# A new perturbative technique for solving integro-partial differential equations

Peter A. Becker[a)]
*Center for Earth Observing and Space Research, Institute for Computational Sciences and Informatics, and Department of Physics and Astronomy, George Mason University, Fairfax, Virginia 22030-4444*

Integro-partial differential equations occur in many contexts in mathematical physics. Typical examples include time-dependent diffusion equations containing a parameter (e.g., the temperature) that depends on integrals of the unknown distribution function. The standard approach to solving the resulting nonlinear partial differential equation involves the use of predictor–corrector algorithms, which often require many iterations to achieve an acceptable level of convergence. In this paper we present an alternative procedure that allows us to separate a family of integro-partial differential equations into two related problems, namely (i) a perturbation equation for the temperature, and (ii) a linear partial differential equation for the distribution function. We demonstrate that the variation of the temperature can be determined by solving the perturbation equation *before* solving for the distribution function. Convergent results for the temperature are obtained by recasting the divergent perturbation expansion as a continued fraction. Once the temperature variation is determined, the self-consistent solution for the distribution function is obtained by solving the remaining, linear partial differential equation using standard techniques. The validity of the approach is confirmed by comparing the (input) continued-fraction temperature profile with the (output) temperature computed by integrating the resulting distribution function.

## I. INTRODUCTION

Many of the time-dependent transport equations encountered in mathematical physics are nonlinear in nature due to the dependence of one or more of the coefficients on integrals of the unknown distribution function. In such cases, the transport equation becomes an integro-partial differential equation such as the Vlasov or Boltzmann equations. Physical applications include a large variety of diffusive and plasma phenomena[1,2] as well as nonlinear wave propagation,[3] the dynamics of self-gravitating mass distributions,[4] and the diffusion in energy space of photons due to Compton scattering.[5,6]

Integro-partial differential equations are usually solved by integrating forward in time from a given initial condition using a predictor–corrector algorithm[7,8] or a global relaxation method.[9,10] The convergence properties of such indirect methods are often difficult to predict in advance, and usually depend rather sensitively on both the governing equation and the nature of the initial conditions. In this paper we develop an alternative procedure that allows us to analyze the time variation of the integral function (in this case the temperature) using a *direct* method based upon the governing integro-partial differential equation. The temperature variation is determined by constructing a perturbation expansion via recursive application of the moment equation obtained by integrating the original nonlinear equation with respect to energy.

The coefficients of the perturbation expansion depend on the initial shape of the distribution

[a)]Electronic mail: pbecker@gmu.edu

function, and therefore the resulting temperature variation represents the self-consistent solution to the problem. With the temperature variation determined in advance, the equation governing the distribution function loses its integrodifferential character, and reduces to a linear partial differential equation which can be solved using a variety of standard techniques. The accuracy of the solution can be verified *a posteriori* by integrating the resulting distribution function over energy to obtain another (output) result for the variation of the temperature integral, which can be compared with the (input) temperature representation constructed using the information contained in the perturbation expansion.

In a certain sense, the method developed here allows us to ''separate'' the original integro-partial differential equation into two problems, the first being the determination of the self-consistent temperature variation and the second the solution of the remaining linear partial differential equation for the distribution function. The perturbation series for the temperature is divergent in general, but we demonstrate that it can be recast as a continued fraction that yields convergent results for a variety of initial distributions. In order to illustrate the technique, we focus here on a specific family of equations which is sufficiently general to admit a variety of interesting behaviors.

## II. GOVERNING EQUATIONS

The primary motivation for this study is the analysis of the scattering of photons and electrons in the hot, tenuous plasma surrounding a compact astrophysical object such as a neutron star or black hole. This process, referred to as time-dependent Comptonization, is thought to be responsible for producing the variable x-ray emission observed from a variety of sources both within and outside our galaxy. The energy of the photons is modified as a result of multiple interactions with electrons, and consequently the photon energy distribution evolves over time. In this situation the photon distribution function $f(x,y)$ is governed by an integro-partial differential transport equation of the general form[5]

$$\frac{\partial f}{\partial y}=\frac{1}{x^i}\,\frac{\partial}{\partial x}\left\{x^i\left[x^j\frac{f}{\theta(y)}+x^k\frac{\partial f}{\partial x}\right]\right\}, \tag{1}$$

where $x$ represents the dimensionless photon energy, $y$ measures the dimensionless time, and $i$, $j$, and $k$ are constants. The function $\theta(y)$ represents the time-varying temperature, defined by the integral expression

$$\theta(y)\equiv\frac{I_\alpha(y)}{I_\alpha(0)}, \tag{2}$$

where $\alpha$ is a constant and the power moments of $f$ are defined by

$$I_n(y)\equiv\int_0^\infty x^n f(x,y)\,dx. \tag{3}$$

Note that $\theta(0)=1$ by virtue of (2), and $\theta(y)>0$ for all $y$ since $f(x,y)$ is non-negative. In the time-dependent Comptonization problem, we have $i=j=k=2$. However, we will develop the formalism for arbitrary values of $i$, $j$, and $k$ in order to emphasize the generality of the mathematical method. For clarity in the discussion, we shall think of the test particles as ''photons'' and the scattering centers as ''electrons,'' although these identifications are arbitrary.

The total number density of the photons $N_r(y)$ is related to the photon distribution function $f(x,y)$ via

$$N_r(y)\equiv\int_0^\infty x^i f(x,y)\,dx, \tag{4}$$

so that $N_r(y)=I_i(y)$. Interpreting $x^{-i}(\partial/\partial x)x^i$ as the divergence operator in energy space, we observe that the transport equation (1) is written in explicit flux-conservation form, and therefore $N_r$ remains constant since (1) contains no sources or sinks of photons. We seek to solve the transport equation for the distribution function $f(x,y)$ subject to the initial condition

$$f(x,0)\equiv f_0(x), \tag{5}$$

where $f_0(x)$ is a known function specified as part of the problem under consideration. The conserved number density is therefore given by $N_r=\int_0^\infty x^i f_0(x)dx$. Equation (1) drives $f(x,y)$ toward the steady state equilibrium solution given by the exponential spectrum

$$f_{\rm eq}(x)\equiv \frac{N_r p}{(p\,\theta_{\rm eq})^{(i+1)/p}\Gamma\left(\dfrac{i+1}{p}\right)}\exp\left(\frac{-x^p}{p\,\theta_{\rm eq}}\right), \quad p\equiv j-k+1, \tag{6}$$

where $\theta_{\rm eq}$ is the asymptotic equilibrium temperature and the normalization has been set so that the number density of $f_{\rm eq}$ is equal to $N_r$. In deriving (6) we have also assumed that $(i+1)/p>0$. Whether or not the solution $f(x,y)$ actually reaches $f_{\rm eq}(x)$ depends upon the rate at which the temperature varies in a given situation.

The underlying process modeled by (1) is a stochastic energization of the test particles due to the random motions of the scattering centers. This interpretation is made clear by using (1) to calculate the Fokker–Planck coefficients which express the rates of change of the mean energy $\langle x\rangle$ and the variance $\sigma^2$ for a monoenergetic distribution. The results obtained are

$$\frac{d\langle x\rangle}{dy}=(i+k)x^{k-1}-\frac{x^j}{\theta(y)}, \quad \frac{d\sigma^2}{dy}=2x^k, \tag{7}$$

which describe, respectively, the ''drifting'' and ''broadening'' of the distribution due to energy space diffusion.[11] In terms of these coefficients, (1) can be recast as the Fokker–Planck equation

$$\frac{\partial F}{\partial y}=-\frac{\partial}{\partial x}\left[F\frac{d\langle x\rangle}{dy}\right]+\frac{\partial^2}{\partial x^2}\left[F\frac{1}{2}\frac{d\sigma^2}{dy}\right], \tag{8}$$

where the photon number spectrum $F(x,y)$ is defined by

$$F(x,y)\equiv x^i f(x,y), \tag{9}$$

so that $N_r=\int_0^\infty F(x,y)dx$. It can be readily verified that (8) is equivalent to (1). Since $i$, $j$, and $k$ are free parameters and $\theta(y)$ is an arbitrary power integral of $F$, we see that the Fokker–Planck coefficients associated with (1) encompass a large variety of microphysical scattering scenarios.

## III. PERTURBATION EXPANSION FOR THE TEMPERATURE

We can obtain a relationship between the power moments by operating on (1) with $\int_0^\infty x^n\,dx$ and integrating by parts twice. This yields

$$\frac{dI_n}{dy}=(n-i)\left[(n+k-1)I_{n+k-2}(y)-\frac{I_{n+j-1}(y)}{\theta(y)}\right], \tag{10}$$

where we have assumed that the power moments $I_n$ exist for all of the required values of $n$. The validity of this assumption depends on the asymptotic behavior of the initial distribution $f_0(x)$. By recursively applying (10), we can express all of the derivatives of any power moment with respect to $y$ as closed functions of the moments.

One interesting consequence is that we can obtain the derivatives of the temperature integral function $\theta(y)$ as functions of the moments $I_n(y)$. Using (2), the zeroth derivative is given by

$$\theta(y)=\frac{I_\alpha(y)}{I_\alpha(0)}. \tag{11}$$

By making a single application of (10), we find that the first derivative can be expressed as

$$\theta^{(1)}(y)=\frac{\alpha-i}{I_\alpha(0)}\left[(\alpha+k-1)I_{\alpha+k-2}(y)-\frac{I_{\alpha+j-1}(y)}{\theta(y)}\right]. \tag{12}$$

Differentiation of (12) with respect to $y$ yields

$$\theta^{(2)}(y)=\frac{\alpha-i}{I_\alpha(0)}\left[\frac{\theta^{(1)}(y)}{\theta^2(y)}I_{\alpha+j-1}(y)-\frac{1}{\theta(y)}\frac{dI_{\alpha+j-1}}{dy}+(\alpha+k-1)\frac{dI_{\alpha+k-2}}{dy}\right]. \tag{13}$$

Using (10) to eliminate the moment derivatives in (13), we obtain

$$\begin{aligned}\theta^{(2)}(y)=\frac{\alpha-i}{I_\alpha(0)}\Bigg\{&(\alpha+k-1)(\alpha+k-i-2)\left[(\alpha+2k-3)I_{\alpha+2k-4}(y)-\frac{I_{\alpha+k+j-3}(y)}{\theta(y)}\right]\\&+\frac{\theta^{(1)}(y)I_{\alpha+j-1}(y)}{\theta^2(y)}+(\alpha+j-i-1)\left[\frac{I_{\alpha+2j-2}(y)}{\theta^2(y)}-(\alpha+j+k-2)\frac{I_{\alpha+j+k-3}(y)}{\theta(y)}\right]\Bigg\}.\end{aligned} \tag{14}$$

Subsequent iterative applications of (10) can be used to derive expressions for the third and higher derivatives of $\theta(y)$.

By evaluating the derivatives sequentially starting with $\theta(y)$, we can calculate as many as desired if the required power moments $I_n(y)$ are known. Although we have no *a priori* means of evaluating the power moments $I_n(y)$ for general values of $y$, we *can* evaluate them for the special case $y=0$ since in this case they correspond to integrals of the known initial distribution $f_0(x)$, i.e.,

$$I_n(0)=\int_0^\infty x^n f(x,0)dx=\int_0^\infty x^n f_0(x)dx. \tag{15}$$

When $y=0$, the general expressions for the derivatives given by (11), (12), and (14) reduce to

$$\theta(0)=1, \tag{16}$$

$$\theta^{(1)}(0)=\frac{\alpha-i}{I_\alpha(0)}[(\alpha+k-1)I_{\alpha+k-2}(0)-I_{\alpha+j-1}(0)], \tag{17}$$

$$\begin{aligned}\theta^{(2)}(0)=\frac{\alpha-i}{I_\alpha(0)}\{&(\alpha+k-1)(\alpha+k-i-2)[(\alpha+2k-3)I_{\alpha+2k-4}(0)-I_{\alpha+k+j-3}(0)]\\&+\theta^{(1)}(0)I_{\alpha+j-1}(0)+(\alpha+j-i-1)[I_{\alpha+2j-2}(0)-(\alpha+j+k-2)I_{\alpha+j+k-3}(0)]\}.\end{aligned} \tag{18}$$

The method can be extended to evaluate the third and higher derivatives of $\theta(y)$ at $y=0$ for a given initial distribution $f_0(x)$.

Let us suppose that for some arbitrary value of $M$, all of the initial derivatives of $\theta(y)$ up to $\theta^{(M)}(0)$ have been determined using the method outlined above. We may then define the associated asymptotic perturbation (Taylor) series for $\theta(y)$ by writing

$$\Phi_N(y)\equiv\sum_{n=0}^{N}\frac{\theta^{(n)}(0)}{n!}y^n, \tag{19}$$

where $N\leqslant M$ indicates the truncation level of the series. Our expectation is that $\Phi_N(y)$ accurately approximates the time variation of the exact solution for $\theta(y)$ within some finite radius of convergence if $N$ is sufficiently large. Based upon the existence of the Taylor series, we conclude that in principle the variation of $\theta(y)$ can be determined *before* solving for the unknown distribution $f(x,y)$. This accomplishes the formal ''separation'' of the original integro-partial differential equation (1) into two problems. The first problem is the determination of the variation of the integral function, which has been achieved (at least formally) by constructing the Taylor series (19). The second problem is the determination of the spectrum $f(x,y)$, which now reduces to the solution of a *linear* partial differential equation since the function $\theta(y)$ appearing in (1) can be approximated using (19). However, the convergence of the power series (19) introduces some potential complications which we address below.

## IV. CONTINUED FRACTION REPRESENTATION

We have established that it is possible to develop a general computational scheme based on (10) that can be used to evaluate the initial derivatives of the self-consistent temperature integral function $\theta(y)$ in terms of the initial moments $I_n(0)$, which are easily computed using (15) once the initial distribution $f_0(x)$ is specified. From knowledge of the initial $\theta$ derivatives we are able to construct the formal Taylor series $\Phi_N(y)$ given by (19). However, a remaining difficulty centers on the convergence of this series. In many cases the radius of convergence turns out to be too small to be of any practical use, and in certain instances it may even vanish. We therefore seek an alternative means for utilizing the asymptotic information contained in the power series coefficients in order to extract a *global* representation of the function.

A global approximation can be constructed by recasting the data in the form of a continued fraction, which is equivalent to the process of Padé approximation.[12] In many instances this is a remarkably successful approach to the global modeling of an unknown function. We define the continued fraction representation using

$$\Psi_N(y)\equiv\cfrac{c_0}{1+\cfrac{c_1y}{1+\cfrac{c_2y}{1+\cdots\cfrac{c_{N-1}y}{1+c_Ny}}}}, \tag{20}$$

where the constants $c_0,\ldots,c_N$ are the continued fraction coefficients, and the truncation level is indicated by the value of $N$. The continued fraction coefficients can be computed using the information contained in the $\theta$ derivatives by employing the standard two-dimensional algorithm described by Baker and Graves-Morris.[13]

To illustrate the flow of the algorithm, we assume that via successive applications of (10) we have evaluated all of the initial derivatives of the temperature integral function $\theta(y)$ up to $\theta^{(M)}(0)$ for some $M$. The algorithm is initialized by setting the zeroth column of the matrix $A_{n,m}$ using

$$A_{0,m}=\frac{\theta^{(m)}(0)}{m!},\quad 0\leqslant m\leqslant M, \tag{21}$$

and the first column of the matrix is calculated subsequently via

$$A_{1,m}=-\frac{A_{0,m+1}}{A_{0,0}},\quad 0\leqslant m\leqslant M-1. \tag{22}$$

The remaining elements in the matrix are obtained recursively using

$$A_{n,m}=\frac{A_{n-2,m+1}}{A_{n-2,0}}-\frac{A_{n-1,m+1}}{A_{n-1,0}}, \quad 0\leqslant m\leqslant M-n, \tag{23}$$

for $2\leqslant n\leqslant M$. The continued fraction coefficients occupy the zeroth row of the matrix, so that

$$c_n=A_{n,0}\,, \quad 0\leqslant n\leqslant M. \tag{24}$$

Note that the coefficient $c_n$ is a function of the initial derivatives $\theta(0)$, $\theta^{(1)}(0),\ldots,\theta^{(n)}(0)$, and therefore the incorporation of higher-order derivatives into the scheme has no effect on the value of $c_n$.

As discussed in Sec. II, the exact solution for the temperature variation $\theta(y)$ must be positive for all $y>0$ because the distribution function $f$ is non-negative. Hence $\theta(y)$ contains no poles in the domain $y>0$, and this must also be true of any acceptable continued fraction approximation. The singularity structure of the continued fraction $\Psi_N(y)$ can be determined by creating an equivalent rational function, and then solving for the zeros of the denominator. The presence of zeros in the domain $y>0$ causes the appearance of extraneous, unphysical poles in $\Psi_N(y)$. These unphysical poles migrate out of the computational domain as the truncation level $N$ increases and $\Psi_N(y)$ approaches the exact solution $\theta(y)$. However, some of the lower order fractions may contain ''defects'' (poles for positive values of $y$), and if so they must be rejected. We consider the convergence properties of the continued fraction sequence $\Psi_0(y), \Psi_1(y),\ldots,\Psi_M(y)$ in Sec. VI, where we treat specific computational examples.

## V. APPLICATION TO TIME-DEPENDENT COMPTONIZATION

In the problem of time-dependent Comptonization that serves as the primary motivation for this study, an intense distribution of radiation is scattered by hot electrons in a tenuous plasma. We suppose that the photons are injected impulsively into the plasma with a specified energy distribution at time $t=0$, and that the subsequent evolution of the distribution occurs as a result of photon–electron scattering. This process naturally leads to the production of time-variable (transient) x-ray emission. During the brightest observed x-ray transients, most of the energy is contained in the radiation field rather than in the gas, which implies that the average photon energy must remain constant even as the shape of the x-ray spectrum changes. In this case the electrons maintain a Maxwellian distribution and act mainly as catalysts in the evolution of the photon distribution, taking energy away from very energetic photons and giving it to low-energy photons until equilibrium is achieved. The electron temperature $T_e$ depends on the shape of the radiation spectrum via an integral expression, and therefore $T_e$ varies as a function of time as we demonstrate below. The equation governing the evolution of the radiation spectrum is consequently integro-partial differential in nature.

The time evolution of the photon spectrum under the influence of Compton scattering in an ionized, homogeneous hydrogen plasma is governed by the Kompaneets equation,[14]

$$\frac{\partial f}{\partial t}=\frac{n_e\sigma_T ckT_e(0)}{m_ec^2}\,\frac{1}{x^2}\,\frac{\partial}{\partial x}\left\{x^4\left[f+\frac{T_e(t)}{T_e(0)}\,\frac{\partial f}{\partial x}\right]\right\}, \tag{25}$$

where $\sigma_T$ is the Thomson scattering cross section, $k$ is Boltzmann's constant, $c$ is the speed of light, and $n_e$ and $m_e$ denote the electron number density and mass, respectively. The dimensionless energy variable $x$ is defined by

$$x\equiv\frac{\epsilon}{kT_e(0)}, \tag{26}$$

where $\epsilon$ is the photon energy and $T_e(0)$ denotes the electron temperature at the beginning of the transient $(t=0)$. The initial spectrum $f_0(x)\equiv f(x,0)$ is assumed to be known, and $f$ is normalized so that the total photon number density is given by

$$N_r(y)=\int_0^\infty x^2 f(x,y)dx=I_2(y), \tag{27}$$

where the power moments of $f$ are defined by $I_n(y)\equiv\int_0^\infty x^n f(x,y)dx$ in accordance with (3). We remind the reader that $N_r$=constant according to the discussion in Sec. II. The associated total photon energy density is given by

$$U_r(y)=\int_0^\infty \epsilon\, x^2 f(x,y)dx=kT_e(0)I_3(y). \tag{28}$$

In this application, the explicit connection between the ''time parameter'' $y$ and the true time $t$ is established by making the definition

$$y(t)\equiv\int_0^t n_e(t)\sigma_T c\,\frac{kT_e(t')}{m_e c^2}\,dt', \tag{29}$$

where we have allowed for the possibility of a time dependence in the electron number density and we have set $y(0)=0$. Using the variables $x$ and $y$ our transport equation (25) can be reexpressed as

$$\frac{\partial f}{\partial y}=\frac{1}{x^2}\frac{\partial}{\partial x}\left\{x^4\left[\frac{f}{\theta(y)}+\frac{\partial f}{\partial x}\right]\right\}, \tag{30}$$

where the temperature function $\theta(y)$ is defined by

$$\theta(y)\equiv\frac{T_e(y)}{T_e(0)}. \tag{31}$$

In order to apply the method developed earlier in the paper to the current problem, we must demonstrate that the definitions for $\theta(y)$ in (2) and (31) are consistent, which can be established by showing that the right-hand side of (31) is proportional to one of the power moments of $f$ in the Comptonization application. We begin by noting that (30) is formally equivalent to our prototype transport equation (1) if $i=j=k=2$. It follows that the power moments $I_n$ must satisfy (10), which in this case reduces to

$$\frac{dI_n}{dy}=(n-2)\left[(n+1)I_n(y)-\frac{I_{n+1}(y)}{\theta(y)}\right]. \tag{32}$$

During a bright x-ray transient, most of the energy density is contained in the radiation field, and therefore the material gas cannot exchange a significant amount of energy with the photons. Consequently the integrated radiation energy density $U_r(y)=kT_e(0)I_3(y)$ should not change as a result of Comptonization, and therefore $I_3(y)$ must remain equal to its initial value $I_3(0)$. Setting $n=3$ in (32), we find that the condition $dI_3/dy=0$ is satisfied if

$$\theta(y)=\frac{I_4(y)}{4I_3(0)}, \tag{33}$$

which implies that $T_e$ equals the inverse-Compton temperature of the radiation spectrum.[6] Since $\theta(0)=1$ by virtue of (31), we can rewrite (33) as

$$\theta(y)=\frac{I_4(y)}{I_4(0)}. \tag{34}$$

Note that the initial spectrum $f_0(x)$ must satisfy the condition $I_4(0)=4I_3(0)$ in order to be consistent with the requirement that $\theta(0)=1$. Equation (34) establishes the integro-partial differential nature of the governing equation (30) for this application. Hence the formal results obtained earlier in the paper can be applied to the problem of time-dependent Comptonization by setting $\alpha=4$ in (2).

The algorithm derived in Sec. III can be used to directly calculate the initial derivatives of $\theta(y)$ at $y=0$. However, in the Comptonization application under consideration here, it is more convenient to derive a differential recurrence relation between the successive moments $I_n$ and $I_{n+1}$ by rearranging (32) to obtain

$$I_{n+1}(y)=\frac{\theta(y)}{2-n}e^{(n+1)(n-2)y}\frac{d}{dy}[e^{-(n+1)(n-2)y}I_n(y)]. \tag{35}$$

Working in terms of the differential operator

$$\mathcal{D}_n\equiv\frac{1}{2-n}e^{(n+1)(n-2)y}\frac{d}{dy}e^{-(n+1)(n-2)y}, \tag{36}$$

we can apply (35) iteratively to find that

$$I_{n+1}(y)=\theta(y)\mathcal{D}_nI_n=\theta(y)\mathcal{D}_n\theta(y)\mathcal{D}_{n-1}\cdots\theta(y)\mathcal{D}_3I_3(y). \tag{37}$$

Since $I_3(y)=I_3(0)=$constant, we can carry out the differentiation to obtain the moments $I_n(y)$ as functions of the derivatives of $\theta(y)$. The first few results are

$$I_4(y)=4\,\theta(y)I_3(0), \tag{38}$$

$$I_5(y)=-\frac{\theta(y)}{2}[4\,\theta^{(1)}(y)-40\theta(y)]I_3(0), \tag{39}$$

$$I_6(y)=-\frac{\theta(y)}{3}[-360\theta^2(y)+76\theta(y)\,\theta^{(1)}(y)-2\,\theta^{(1)^2}(y)-2\,\theta(y)\,\theta^{(2)}(y)]I_3(0). \tag{40}$$

Similar results can be obtained for the higher moments. Note that $I_n(y)$ is a function of all of the derivatives of $\theta(y)$ up to $\theta^{(n-4)}(y)$. We have developed a computer algorithm based on (37) that efficiently derives expressions for the moments $I_n(y)$ in terms of the derivatives of $\theta(y)$. Since the values of the initial moments $I_n(0)$ are easily calculated using (15) for any initial spectrum $f_0(x)$, these expressions allow us to compute the corresponding initial derivatives $\theta^{(n)}(0)$ sequentially, beginning with the zeroth derivative which is set by (38). The development of the general expressions giving the initial moments as function of the derivatives $\theta^{(n)}(0)$ is the costliest part of the solution procedure. However, once these expressions have been established, they can be used to evaluate the initial derivatives for a variety of different initial spectra at very low cost. The initial derivatives are subsequently used to calculate the continued fraction coefficients using the algorithm discussed in Sec. IV, and the sequence of continued fraction approximations $\Psi_N(y)$ is evaluated using (20). This procedure forms the basis for the computational results presented in Sec. VI.

## VI. COMPUTATIONAL EXAMPLES

In this section we apply our method to obtain quantitative results for the problem of time-dependent astrophysical Comptonization. Our computational procedure is as follows. First we calculate the initial $\theta$ derivatives and the associated continued fraction coefficients using the algorithms discussed in Secs. IV and V. In the second step we use this information to analyze the convergence properties of the sequence of continued fractions $\Psi_N(y)$ and compare them with the corresponding Taylor series $\Phi_N(y)$. In the third step we use a high-order continued fraction to approximate the temperature integral function $\theta(y)$ in the transport equation (30) and then we numerically solve the transport equation for the photon energy distribution $f(x,y)$. Finally, in the fourth step we compare the input temperature function $\Psi_N(y)$ with the result obtained for $\theta(y)$ by integrating the numerical solution for $f(x,y)$. Agreement between these two representations of the temperature confirms the accuracy of the method.

At the beginning of the x-ray transient, the radiation spectrum is given by the initial energy distribution $f_0(x)$ introduced in (5). According to the analysis presented in the preceding sections, knowledge of $f_0(x)$ is sufficient to determine the time variation of the temperature integral function $\theta(y)$. In typical astrophysical situations, the initial energy distributions of greatest interest are the optically thin electron–proton bremsstrahlung (free–free) spectrum

$$f_0(x) = x^{-3} e^{-x/4}, \tag{41}$$

and the monoenergetic spectrum

$$f_0(x) = N_0 \, x_0^{-2} \delta(x - x_0), \tag{42}$$

where $N_0$ is the number density of the photons. In the latter case, we must set $x_0 = 4$ in order to satisfy the condition $\theta(0) = I_4(0)/4I_3(0) = 1$ as required by (31) and (33).

The transport equation (30) drives the distribution $f(x,y)$ toward the steady state equilibrium solution given by the Wien spectrum

$$f_{\rm eq}(x) \equiv \frac{N_r}{2\,\theta_{\rm eq}^3} e^{-x/\theta_{\rm eq}} \ , \tag{43}$$

which has been obtained by setting $i = j = k = 2$ in (6). Although $f$ is always driven toward Wien form at all values of $y$, it may or may not reach equilibrium depending on the shape of the initial spectrum and the corresponding rate at which the temperature varies. This question can be resolved by calculating the asymptotic temperature $\theta_{\rm eq}$ using information contained in the initial spectrum $f_0(x)$. In equilibrium, the dimensionless mean photon energy is given by

$$\overline{x} = \frac{\int_0^\infty x^3 f_{\rm eq}(x)\,dx}{\int_0^\infty x^2 f_{\rm eq}(x)\,dx} = 3\,\theta_{\rm eq}\,, \tag{44}$$

where we have substituted for $f_{\rm eq}(x)$ using (43) to obtain the final result. Conservation of the photon number and energy densities implies that the value of $\overline{x}$ must be conserved, so that we can also write

$$\overline{x} = \frac{\int_0^\infty x^3 f_0(x)\,dx}{\int_0^\infty x^2 f_0(x)\,dx} = \frac{I_3(0)}{I_2(0)}. \tag{45}$$

Equations (44) and (45) can be combined to calculate the asymptotic temperature $\theta_{\rm eq}$ for any initial spectrum $f_0(x)$. In the case of a monoenergetic initial spectrum with $x_0 = 4$, we obtain $\overline{x} = 4$ and therefore $\theta_{\rm eq} = 4/3$. Conversely, in the case of a bremsstrahlung initial spectrum, we obtain $\overline{x} = 0$ because the number density of photons $N_r = I_2(0) = \int_0^\infty x^2 f_0(x)\,dx$ is formally infinite. This implies that in the bremsstrahlung case, $\theta_{\rm eq} = 0$, and therefore no meaningful steady state exists according to (43). It is important to emphasize that our calculation of $\theta_{\rm eq}$ has utilized an energy conservation principle that may not be available in all physical applications of the general transport equation (1).

TABLE I. Numerical results.

| | Monoenergetic spectrum | | Bremsstrahlung spectrum | |
|---|---|---|---|---|
| $n$ | $\theta^{(n)}(0)$ | $c_n$ | $\theta^{(n)}(0)$ | $c_n$ |
| 0 | $1.00\times10^{0}$ | 1.00 | $1.00\times10^{0}$ | 1.00 |
| 1 | $2.00\times10^{0}$ | −2.00 | $-6.00\times10^{0}$ | 6.00 |
| 2 | $-1.20\times10^{1}$ | 5.00 | $1.32\times10^{2}$ | 5.00 |
| 3 | $8.00\times10^{0}$ | −1.67 | $-6.36\times10^{3}$ | 11.13 |
| 4 | $1.87\times10^{3}$ | 3.59 | $5.29\times10^{5}$ | 8.99 |
| 5 | $-2.99\times10^{4}$ | −2.56 | $-6.68\times10^{7}$ | 15.43 |
| 6 | $-6.85\times10^{5}$ | 4.69 | $1.18\times10^{10}$ | 12.28 |
| 7 | $4.07\times10^{7}$ | −3.56 | $-2.76\times10^{12}$ | 19.62 |
| 8 | $3.65\times10^{8}$ | 4.13 | $8.24\times10^{14}$ | 15.72 |
| 9 | $-9.25\times10^{10}$ | −3.33 | $-3.04\times10^{17}$ | 23.44 |
| 10 | $3.42\times10^{11}$ | 5.03 | $1.36\times10^{20}$ | 19.48 |
| 11 | $3.56\times10^{14}$ | −4.77 | $-7.19\times10^{22}$ | 26.87 |
| 12 | $-6.36\times10^{15}$ | 4.53 | $4.47\times10^{25}$ | 23.51 |
| 13 | $-2.20\times10^{18}$ | −4.30 | $-3.22\times10^{28}$ | 30.09 |
| 14 | $7.14\times10^{19}$ | 5.17 | $2.66\times10^{31}$ | 27.62 |
| 15 | $2.10\times10^{22}$ | −5.67 | $-2.49\times10^{34}$ | 33.35 |
| 16 | $-9.31\times10^{23}$ | 4.97 | $2.64\times10^{37}$ | 31.59 |
| 17 | $-2.96\times10^{26}$ | −5.33 | $-3.13\times10^{40}$ | 36.81 |
| 18 | $1.48\times10^{28}$ | 5.20 | $4.13\times10^{43}$ | 35.34 |
| 19 | $5.90\times10^{30}$ | −6.36 | $-6.04\times10^{46}$ | 40.52 |
| 20 | $-2.70\times10^{32}$ | 6.06 | $9.73\times10^{49}$ | 38.88 |
| 21 | $-1.60\times10^{35}$ | −10.63 | $-1.72\times10^{53}$ | 44.40 |
| 22 | $4.98\times10^{36}$ | −7.69 | $3.32\times10^{56}$ | 42.26 |
| 23 | $5.87\times10^{39}$ | −32.83 | $-6.99\times10^{59}$ | 48.45 |
| 24 | $9.03\times10^{40}$ | 23.42 | $1.59\times10^{63}$ | 45.34 |

The values of the initial $\theta$ derivatives and the associated continued fraction coefficients obtained for the bremsstrahlung and monoenergetic initial spectra are presented in Table I. Note the rapid divergence of the derivatives in each case, implying a limited radius of convergence for the Taylor series $\Phi_N(y)$ given by (19). We compare the convergence properties of the Taylor series and the continued fractions for the bremsstrahlung and monoenergetic initial spectra below.

### A. Monoenergetic initial spectrum

The first and second columns of Table I contain, respectively, the results obtained in the monoenergetic case for the initial derivatives $\theta(0)$, $\theta^{(1)}(0),\ldots,\theta^{(M)}(0)$ and the continued fraction coefficients $c_0, c_1,\ldots,c_M$ for $M=24$. This choice for $M$ is arbitrary and of no particular significance. The corresponding sequence of truncated Taylor series $\Phi_N(y)$ is plotted in Fig. 1. The radius of convergence of the Taylor series is limited to $y\lesssim0.15$ even for large $N$, which reflects the rapid increase in the absolute value of the derivatives of $\theta(y)$ at $y=0$. By contrast, the continued fraction coefficients obtained in the monoenergetic case grow much more slowly in absolute value. The continued fraction sequence $\Psi_N(y)$ contains defects (extraneous poles) for some odd values of $N$, but $\Psi_N(y)$ converges (albeit nonuniformly) for even values of $N$, as can be seen in Fig. 2, where we plot the sequence of continued fractions for $N=18$, 20, 22, 24. It is clear from Fig. 2 that the continued fractions converge even for values of $y$ far outside the radius of convergence of the Taylor series. In our search for an accurate approximation to the exact solution for $\theta(y)$, we select the continued fraction $\Psi_N(y)$ that most closely approaches the correct asymptotic value $\theta_{\rm eq}=4/3$ for large $y$. According to Fig. 2, the best agreement is obtained by using the highest order fraction analyzed in this example, which is $\Psi_{24}(y)$.

With $\theta(y)$ approximated using the continued fraction $\Psi_{24}(y)$, the transport equation (30) reduces to a linear, second-order partial differential equation. We solve this equation using the IMSL routine DMOLCH over the range $0 \leqslant y \leqslant 2$ and $0 \leqslant x \leqslant 50$. The monoenergetic initial condition (42) imposed at $y=0$ is approximated using a Gaussian distribution with mean $\bar{x}=4$, variance $\sigma^2=0.01$, and photon number density $N_0=1$. The result obtained for the photon energy spectrum $G(x,y) \equiv x^3 f(x,y)$ is plotted in Fig. 3. It is convenient to plot $G$ rather than $f$ because energy conservation implies that $I_3$=constant, and therefore $\int_0^\infty G\,dx$ is independent of $y$. In the monoenergetic case the initial condition (42) yields $I_3=4$. As $y$ increases from zero, the distribution evolves away from the monoenergetic initial form and is well described by the equilibrium Wien distribution (43) for $y \gtrsim 1$.

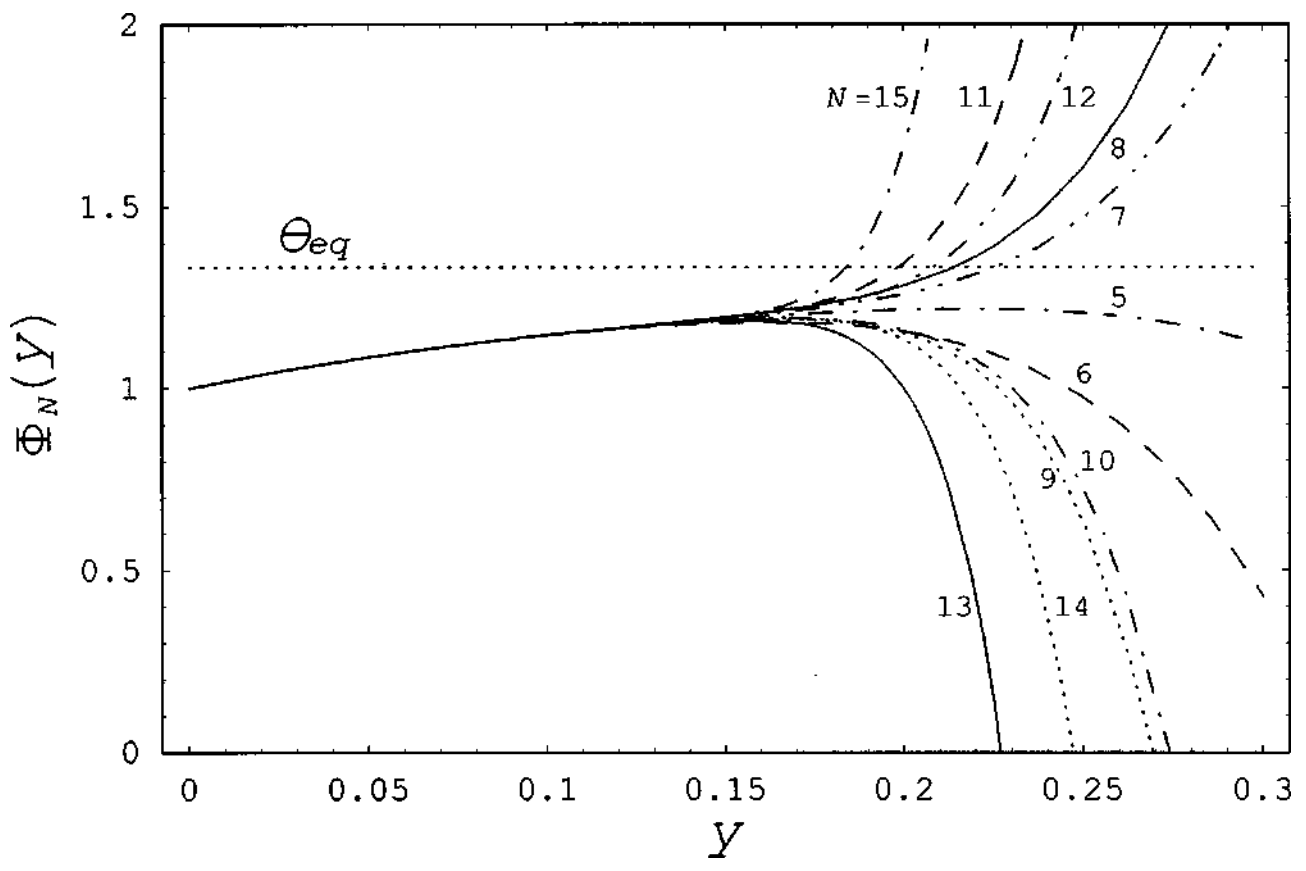

FIG. 1. Sequence of truncated Taylor series $\Phi_N(y)$ given by (19) obtained in the case of Comptonization of a monoenergetic initial spectrum (42) plotted as a function of $y$, with the truncation level $N$ indicated for each curve. Note that $\Phi_N(y)$ diverges for $y \gtrsim 0.15$ even for large $N$, due to the rapid growth of the magnitude of the initial $\theta$ derivatives, as can be seen in Table I. The asymptotic equilibrium temperature $\theta_{eq}=4/3$ is denoted by the dotted horizontal line.

## B. Bremsstrahlung initial spectrum

The third and fourth columns of Table I contain, respectively, the results obtained in the bremsstrahlung case for the initial derivatives $\theta(0)$, $\theta^{(1)},\ldots,\theta^{(M)}(0)$ and the continued-fraction coefficients $c_0, c_1, \ldots, c_M$ for $M=24$. The corresponding sequence of truncated Taylor series $\Phi_N(y)$ is plotted in Fig. 4, and the sequence of continued fractions $\Psi_N(y)$ is plotted in Fig. 5. Note that the radius of convergence of the Taylor series actually *decreases* with increasing truncation level $N$. Conversely, the sequence of continued fractions displays a pattern of uniform convergence with increasing $N$. The uniform convergence is a consequence of the fact that the computed continued fraction coefficients are all positive in this case, leading us to conjecture that the exact solution $\theta(y)$ is a Stieltjes function when the initial spectrum corresponds to optically thin bremsstrahlung.[12]

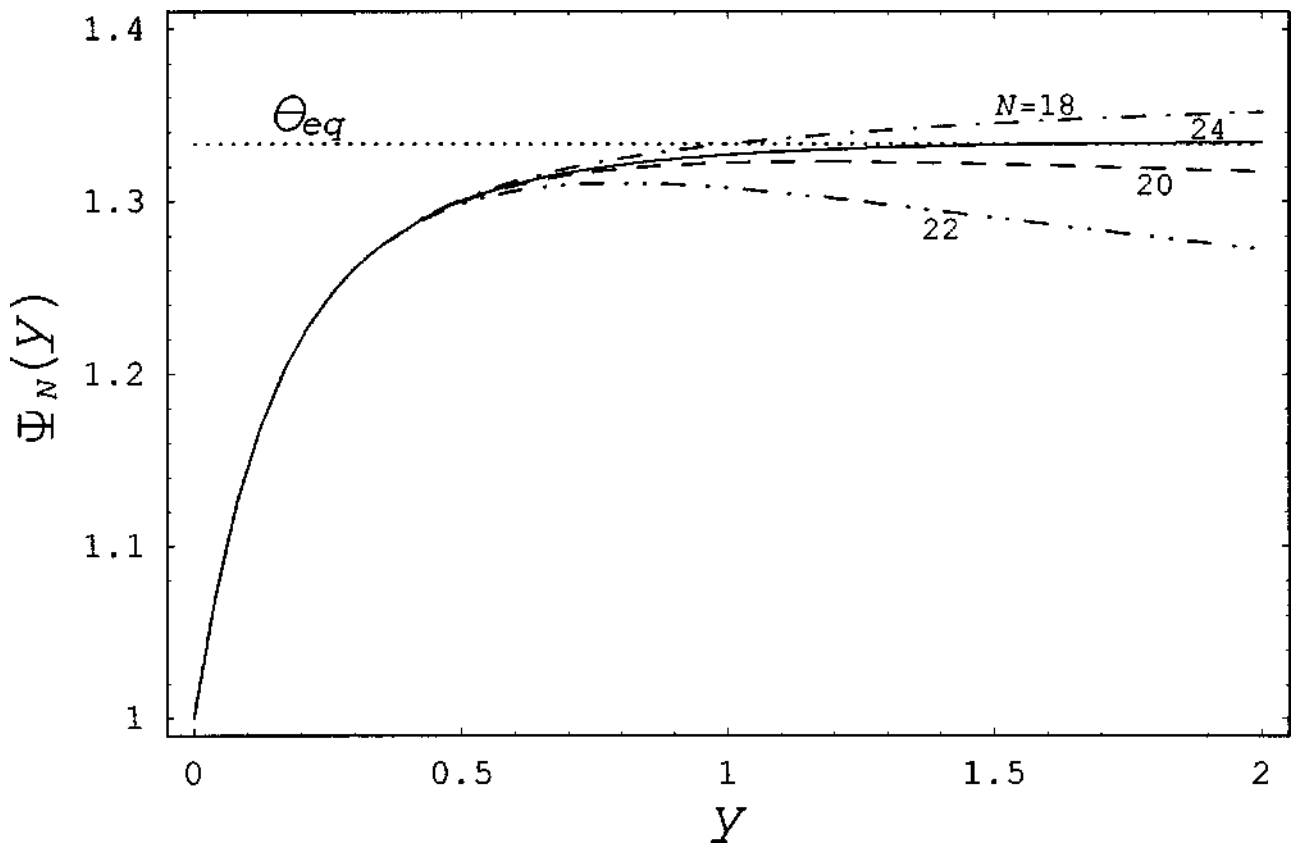

FIG. 2. Sequence of continued fractions $\Psi_N(y)$ given by (20) obtained in the case of Comptonization of a monoenergetic initial spectrum (42) plotted as a function of $y$, with the truncation level $N$ indicated for each curve. We have plotted the results only for even values of $N$ because in this example $\Psi_N(y)$ contains defects (extraneous poles) for some odd values of $N$. Note that $\Psi_N(y)$ converges (though nonuniformly) for even values of $N$ throughout the entire computational domain, which extends well beyond the radius of convergence of the Taylor series depicted in Fig. 1. We require that acceptable approximations approach $\theta_{eq}$ as $y \rightarrow \infty$.

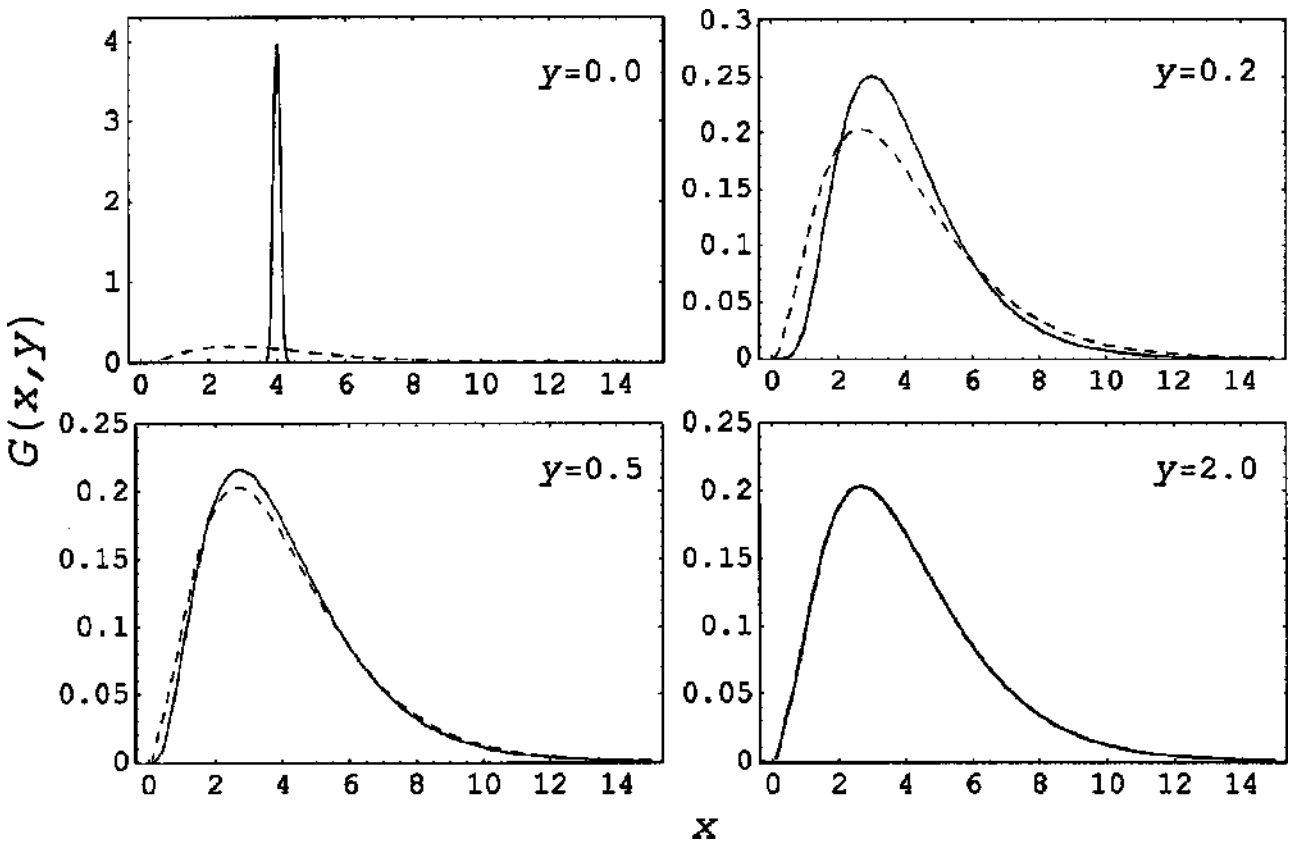

FIG. 3. Numerical result for the photon distribution $G(x,y)\equiv x^3 f(x,y)$ plotted as a function of $x$ for the indicated values of $y$. In this case the initial spectrum imposed at $y=0$ corresponds to a Gaussian distribution with mean $\overline{x}=4$, variance $\sigma^2=0.01$, and total photon number density $N_0=1$, which approximates the monoenergetic initial spectrum (42). The solid lines denote $G$ and the dashed lines represent the asymptotic equilibrium Wien spectrum given by (43). The photon spectrum is essentially given by the Wien form for $y \gtrsim 1$, and the two distributions are indistinguishable for $y=2$. The area under the curves $\int_0^\infty G\,dx=4$ due to energy conservation.

The pattern of uniform convergence of $\Psi_N(y)$ in the bremsstrahlung case suggests that we can obtain a reasonable approximation for the exact solution $\theta(y)$ using the continued fraction $\Psi_{24}(y)$. We impose the bremsstrahlung initial condition (41) at $y=0$ and solve the transport equation (30) using DMOLCH over the range $0\leqslant y\leqslant 2$ and $0\leqslant x\leqslant 50$. The result obtained for the function $G(x,y)\equiv x^3 f(x,y)$ is plotted in Fig. 6. The initial condition (41) combined with energy conservation implies that the area under the curve $\int_0^\infty G\,dx=4$ for all $y$. As $y$ increases from zero, the distribution evolves away from its initial bremsstrahlung form and attempts to approach the equilibrium Wien distribution (43). However, in this case the photon number density $N_r \rightarrow \infty$, and therefore the asymptotic value for the temperature $\theta_{\rm eq} = 0$ as discussed earlier. The divergence of the number density is due to the presence of an infinite number of zero-energy photons in the initial spectrum. Since these photons cannot all be upscattered to higher energies without violating overall energy conservation, the ''equilibrium'' solution for the distribution function given by (43) reduces to a pulse centered on zero energy, which is not a meaningful steady state solution. Hence equilibrium cannot be achieved in the bremsstrahlung case, in contrast to the result obtained when the initial spectrum contains a finite number density of photons, as in the monoenergetic example treated above. We analyze the self-consistency of the numerical solutions obtained in the monoenergetic and bremsstrahlung cases in Sec. VII.

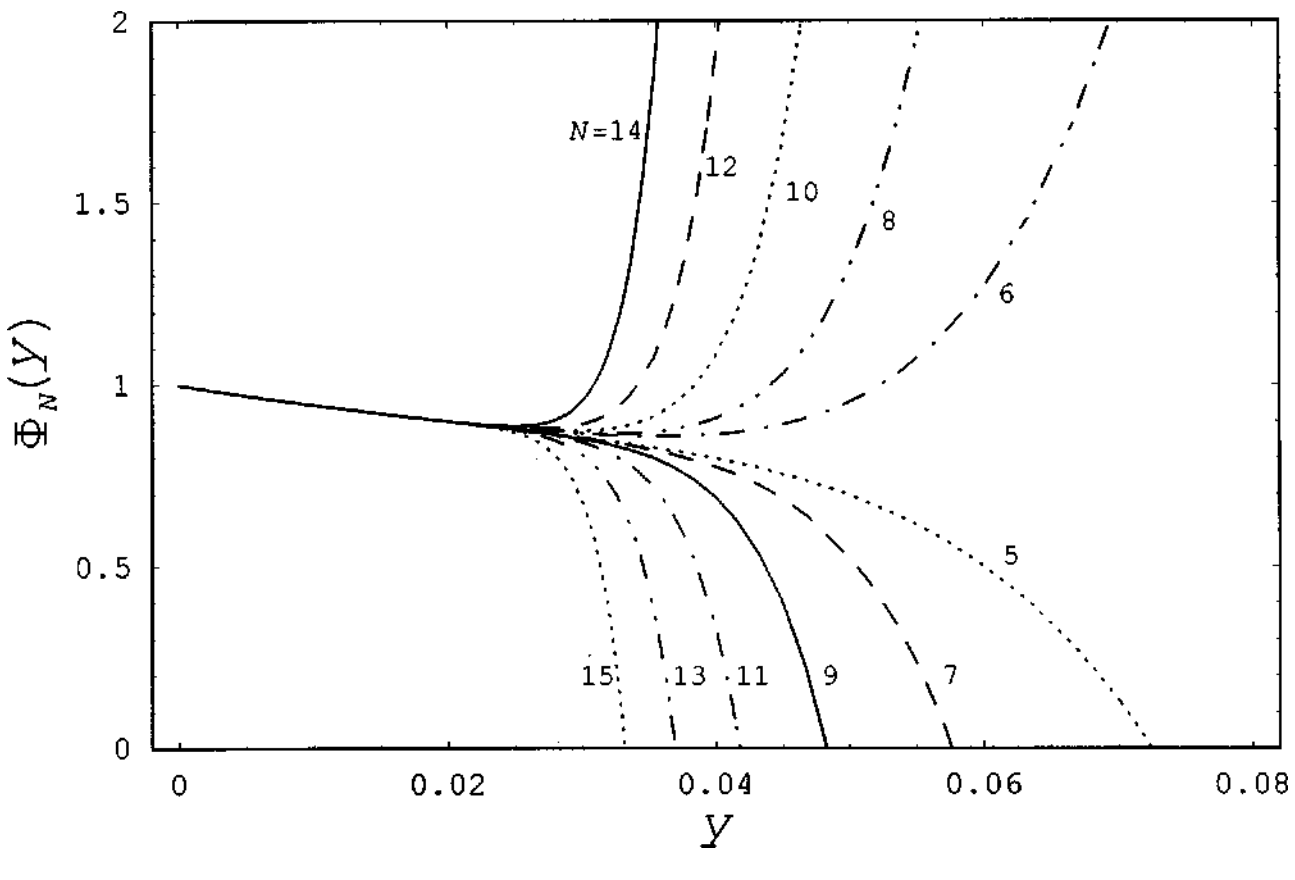

FIG. 4. Sequence of truncated Taylor series $\Phi_N(y)$ obtained in the case of Comptonization of a bremsstrahlung initial spectrum (41) plotted as a function of $y$, with the truncation level $N$ indicated for each curve. Note that the radius of convergence decreases with increasing $N$ due to the rapid growth of the initial $\theta$ derivatives (see Table I). In this case the asymptotic equilibrium temperature $\theta_{eq}=0$ due to the presence of an infinite number of zero-energy photons in the initial spectrum.

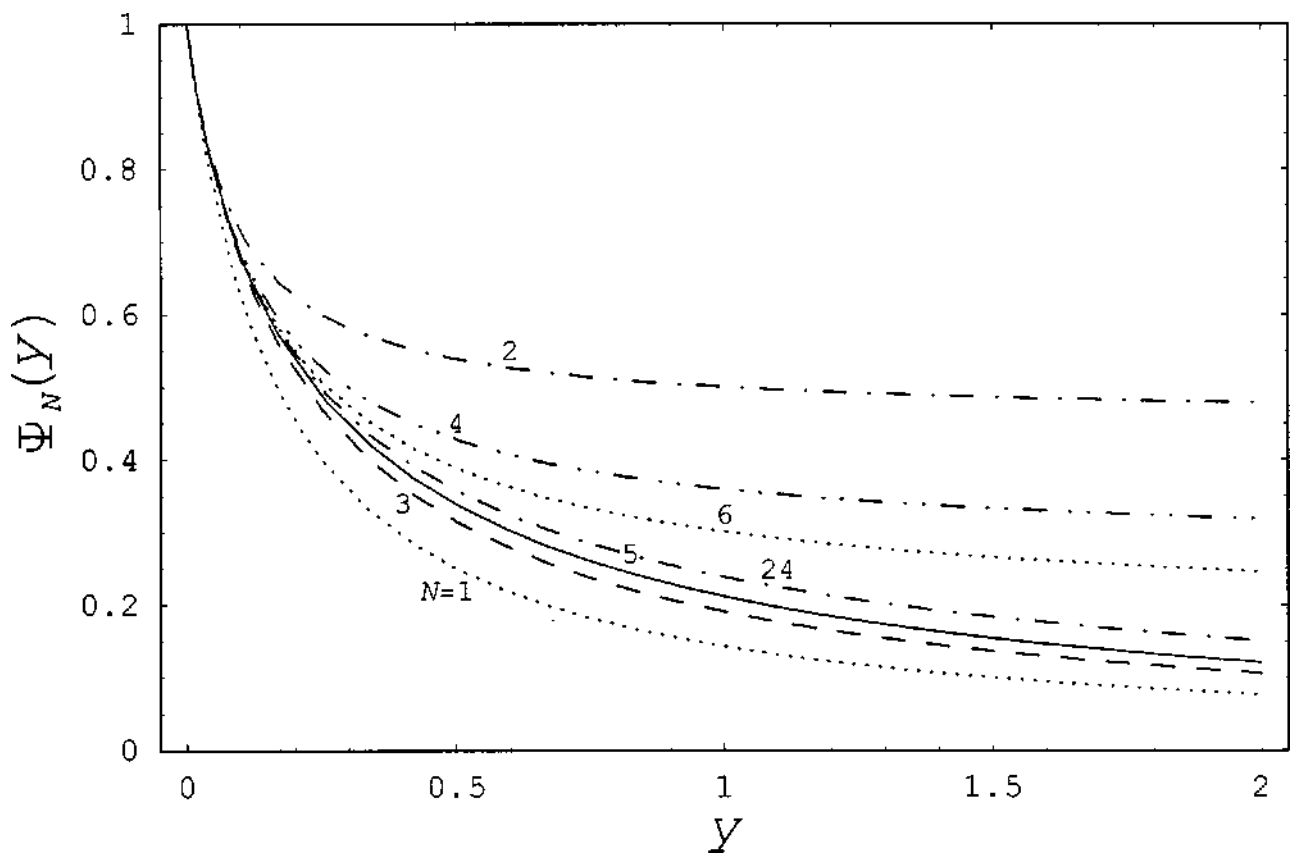

FIG. 5. Sequence of continued fractions $\Psi_N(y)$ obtained in the case of a bremsstrahlung initial spectrum (41) plotted as a function of $y$, with the truncation level $N$ indicated for each curve. Note the pattern of uniform convergence for $N = 1,2,3,4,5,6$, which suggests that in this case the exact solution $\theta(y)$ is probably a Stieltjes function. For $N > 6$, the results become strongly clustered around the $N = 24$ curve.

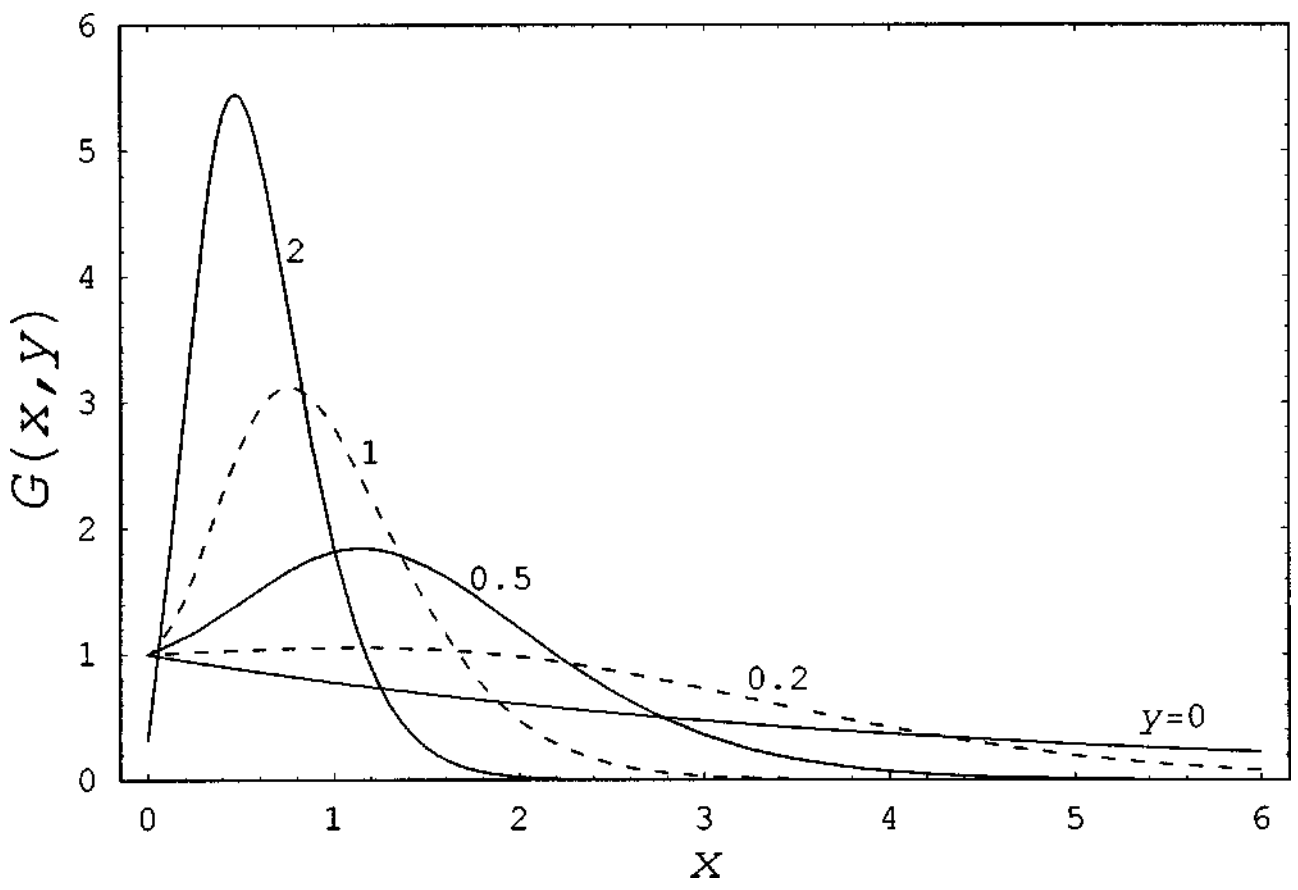

FIG. 6. Numerical results for the photon distribution $G(x,y) \equiv x^3 f(x,y)$ plotted as a function of $x$ for the indicated values of $y$. In this case the initial condition is the optically thin bremsstrahlung spectrum (41). The distribution function approaches a pulse centered at zero energy as $y$ increases due to the infinite number of zero-energy photons in the initial spectrum. The area under the curves $\int_0^\infty G\,dx = 4$ as a consequence of energy conservation.

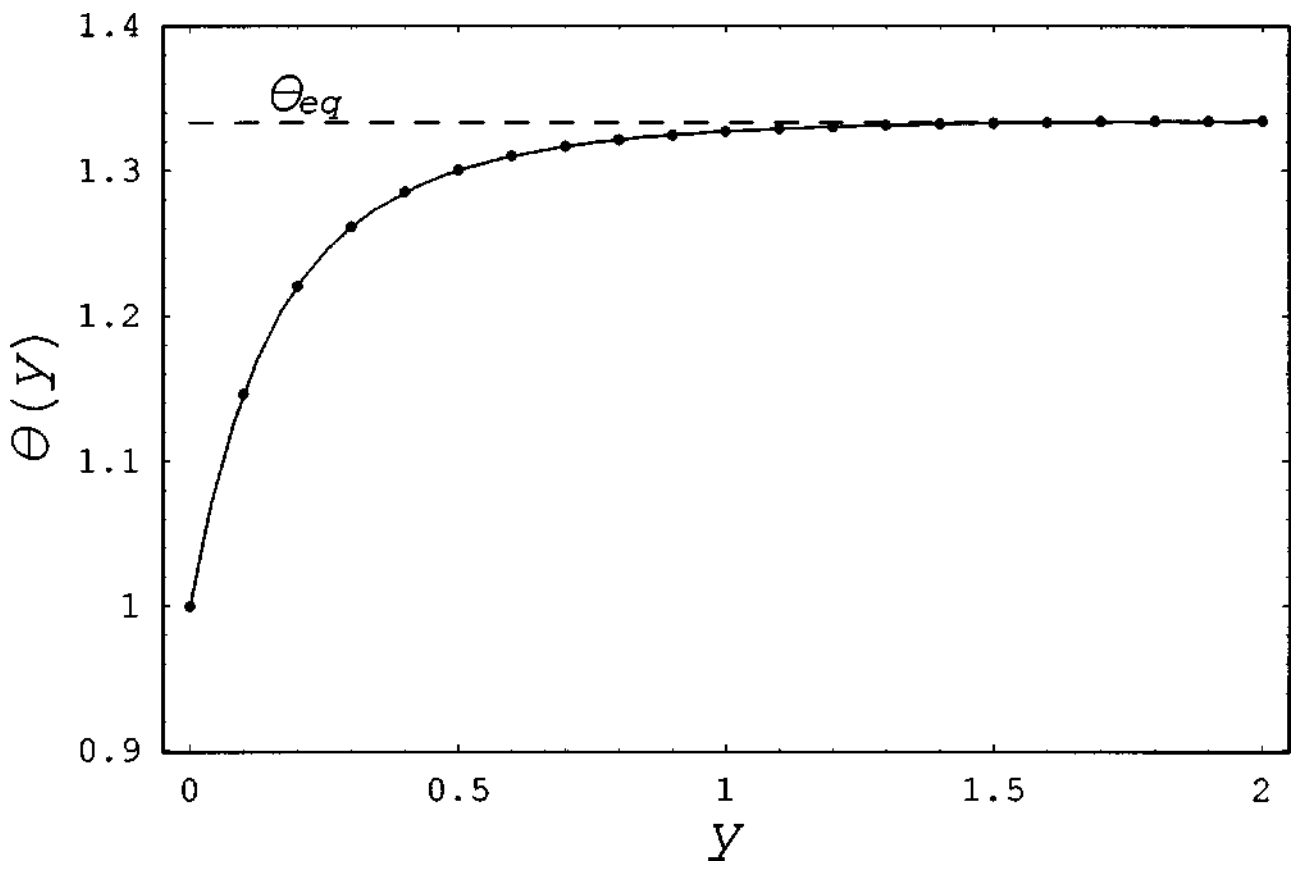

FIG. 7. Results obtained for the (input) continued fraction approximation $\Psi_{24}(y)$ (solid line) are compared with the (output) temperature calculated using the integral definition $\theta(y) \equiv I_4(y)/I_4(0)$ (closed dots) for the Comptonization of a monoenergetic initial spectrum. The dashed horizontal line indicates the asymptotic value of the temperature. The agreement between the two results confirms the self-consistency of the solution obtained for the spectrum $f(x,y)$, which validates the mathematical approach.

## VII. EVALUATION OF THE METHOD

Our approach to the solution of the integro-partial differential transport equation (30) has been to approximate the exact solution for the temperature function $\theta(y)$ by using the available derivatives to construct the highest-order continued fraction that is consistent with the known asymptotic behavior of the temperature. With $\theta(y)$ approximated in this way, we have solved for the photon distribution $f(x,y)$ using a standard computer algorithm commonly available in the IMSL library. One may well ask whether there is any guarantee that the numerical solution so obtained is actually the correct physical solution. Interestingly, for the problem treated here there *is* a method that can be used to *guarantee* both the accuracy and the uniqueness of the solution. This is accomplished by integrating the numerically obtained distribution $f(x,y)$ to calculate the corresponding temperature distribution *a posteriori* using (2) and (3). In the case of astrophysical Comptonization that serves as our sample application in this paper, the corresponding expression is $\theta(y) = I_4(y)/I_4(0)$ as given by (34). A comparison between the result for $\theta(y)$ obtained in this manner and the continued fraction approximation $\Psi_N(y)$ used in the solution of the transport equation serves as the acid test of the entire mathematical and computational approach presented here. In Figs. 7 and 8 we perform this comparison for the bremsstrahlung and monoenergetic initial spectra, respectively. The agreement is clearly excellent, verifying the validity of the overall approach.

## VIII. CONCLUSION

The technique developed here provides a powerful tool for determining the time dependence of the integral function in an integro-partial differential equation *before* solving for the unknown distribution. The numerical results we have obtained for the variation of the self-consistent temperature function $\theta(y)$ in the case of astrophysical Comptonization suggest that the method has acceptable accuracy and reliable convergence properties. The first step in the procedure is the determination of the initial derivatives of $\theta(y)$ at $y = 0$ using the algorithm described in Sec. III, which is based on the differential equation (10) governing the power moments $I_n(y) = \int_0^\infty x^n f(x,y) dx$. The initial derivatives $\theta^{(n)}(0)$ are then used to calculate the continued fraction coefficients $c_n$ appearing in the representation of the continued fraction $\Psi_N(y)$. The convergence of the continued fraction sequence is then analyzed and the highest-order fraction that can be constructed using the available set of coefficients is used to approximate the exact solution $\theta(y)$. Next, the transport equation (1) is solved numerically, using the continued fraction approximation as input to calculate the integral term. Finally, the self-consistency of the solution is evaluated by comparing the input continued fraction approximation $\Psi_N(y)$ with the output temperature calculated using the integral definition $\theta(y) \equiv I_\alpha(y)/I_\alpha(0)$. The level of agreement between $\Psi_N(y)$ and $\theta(y)$ provides a measure of the overall accuracy of the solution procedure.

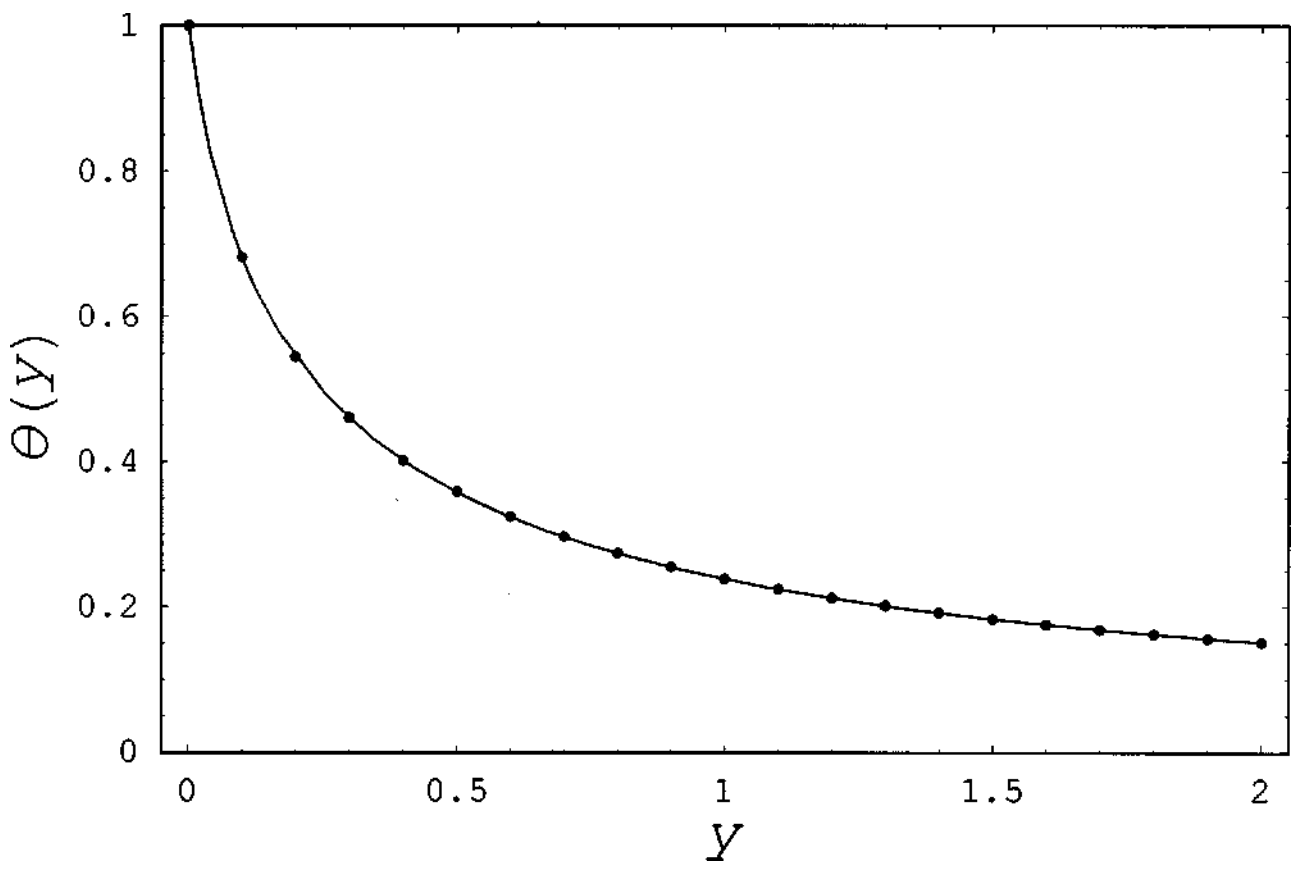

FIG. 8. Results obtained for the (input) continued fraction approximation $\Psi_{24}(y)$ (solid line) are compared with the (output) temperature calculated using the integral definition $\theta(y) \equiv I_4(y)/I_4(0)$ (closed dots) for the Comptonization of a bremsstrahlung initial spectrum. The agreement between the two results confirms the self-consistency of the solution.

It is interesting to contrast the behavior of the continued fraction sequence $\Psi_N(y)$ with that of the associated Taylor series $\Phi_N(y)$. Since the self-consistent solution for the temperature integral function $\theta(y)$ is positive for all real $y>0$, the limited radius of convergence of the Taylor series must reflect the presence of poles in $\theta(y)$ somewhere else in the complex plane. The existence of these poles essentially dooms any attempt to construct a useful perturbation series for the temperature. Conversely, the success of the continued fraction representation stems from its ability to produce a convergent *global* function that shares the same singularity structure as the exact solution $\theta(y)$. In general, $\Psi_N(y)$ converges toward the exact solution as the truncation level $N$ increases, although the pattern of convergence can vary significantly from problem to problem. The cases for which $\theta(y)$ is a Stieltjes function are of particular importance because in these situations the convergence of the continued fraction sequence is uniform.[12]

In our computational examples, which focus on astrophysical Comptonization, we are able to utilize energy conservation to derive the asymptotic value of the temperature $\theta_{\text{eq}}$ in the limit $y \to \infty$. This type of asymptotic information may not be available in every physical application of the general transport equation (1), but if $\theta_{\text{eq}}$ can be calculated, we also have the option of including this information directly into the continued fraction using a two-point algorithm such as those given by Becker[15] and by Baker and Graves-Morris.[13] When this information is incorporated into the continued fraction, $\Psi_N(y)$ automatically approaches $\theta_{\text{eq}}$ as $y \to \infty$. Note that in order to construct the two-point continued fraction we must first transform the time variable from the infinite domain $0 \leqslant y < \infty$ to an equivalent finite domain.

The method presented here bears some relation to techniques for solving partial differential equations proposed by Jumarie[16] and by Bender, Boettcher, and Milton.[17] In Jumarie's approach, a linear Fokker–Planck equation is used to generate moment equations similar to ours, which are solved using a maximum entropy principle to obtain the distribution function. Conversely, Bender, Boettcher, and Milton determine the distribution function governed by a nonlinear partial differential equation by employing a perturbation expansion followed by Padé summation, which resembles our approach to modeling the integral function $\theta(y)$. Although the procedures developed by these authors incorporate certain elements of the technique presented here, a crucial distinction is that their methods are not applicable to integro-partial differential equations.

In conclusion we point out that the procedure developed in this paper can be used as the basis for a new solution technique, or it can be incorporated into existing predictor–corrector or global iteration algorithms as a means of generating a trial solution for the integral function, which is subsequently improved upon using the standard methods. Our focus here has been upon the details of the method, and therefore we have not presented any comparisons between the efficiencies of the various algorithms available for solving integro-partial differential equations. Nonetheless, it is reasonable to expect that the method proposed here is likely to be quite efficient because the temperature variation is determined in advance of solving for the distribution function. In our computational examples we have treated the problem of astrophysical Comptonization, which is governed by the transport equation (30). However, we emphasize that the method is applicable to any equation of the form represented by (1), and that it can potentially be generalized to treat other transport equations, including those containing inhomogeneous, time-dependent source terms.

## ACKNOWLEDGMENTS

The author is grateful to Prasad Subramanian for assistance in carrying out the computational work, and also acknowledges several useful conversations with Ben Crain and Todd Peterson.